\documentclass[journal,twoside,web]{ieeecolor}

\usepackage{generic,lcsys}
\makeatletter
\let\NAT@parse\undefined
\makeatother

\usepackage{amsmath,amssymb,amsthm}
\usepackage{graphicx}
\usepackage{eso-pic}
\usepackage{algorithm,algorithmic}
\usepackage{aliascnt}
\usepackage[sort]{cite}
\usepackage[hidelinks]{hyperref}
\usepackage[capitalise]{cleveref}
\usepackage{etoolbox}

\makeatletter
\patchcmd{\@makecaption}{\setbox\@tempboxa\hbox}{\sbox\@tempboxa}{}%
  {\PackageError{clqr}{Cannot patch the caption measurement box}{}}
\patchcmd{\@makecaption}{\setbox\@tempboxa\hbox}{\sbox\@tempboxa}{}%
  {\PackageError{clqr}{Cannot patch the caption label box}{}}
\makeatother

\newtheorem{theorem}{Theorem}
\newaliascnt{lemma}{theorem}
\newtheorem{lemma}[lemma]{Lemma}
\aliascntresetthe{lemma}
\newaliascnt{proposition}{theorem}
\newtheorem{proposition}[proposition]{Proposition}
\aliascntresetthe{proposition}
\newaliascnt{corollary}{theorem}
\newtheorem{corollary}[corollary]{Corollary}
\aliascntresetthe{corollary}
\newtheorem{definition}{Definition}

\DeclareMathOperator{\rank}{rank}
\DeclareMathOperator{\range}{range}
\DeclareMathOperator{\nullsp}{null}
\newcommand{\R}{\mathbb{R}}
\newcommand{\trans}{\mathsf{T}}

\title{An Efficient Parallel Solver for Linear-Equality-Constrained LQR Problems}

\author{João Sousa-Pinto and Dominique Orban%
\thanks{João Sousa-Pinto (e-mail: \href{mailto:joaospinto@gmail.com}{joaospinto@gmail.com})
carried out this work while at Tesla. Dominique Orban
(e-mail: \href{mailto:dominique.orban@gerad.ca}{dominique.orban@gerad.ca})
is with GERAD and the Department of Mathematics and Industrial Engineering
at Polytechnique Montréal.}\\}

\begin{document}
\maketitle
\thispagestyle{empty}

% IEEE notice for posting the accepted manuscript on arXiv.
\AddToShipoutPictureFG*{%
  \AtPageLowerLeft{%
    \put(\LenToUnit{\dimexpr1in+\hoffset+\oddsidemargin\relax},\LenToUnit{12pt}){%
      \parbox[b]{\textwidth}{%
        \normalfont\fontsize{7}{8}\selectfont\color{black}\raggedright
        \textcopyright~2026 IEEE. Personal use of this material is permitted.
        Permission from IEEE must be obtained for all other uses, in any current
        or future media, including reprinting/republishing this material for
        advertising or promotional purposes, creating new collective works,
        for resale or redistribution to servers or lists, or reuse of any
        copyrighted component of this work in other works.}}}}

\begin{abstract}
We present a parallel solver for finite-horizon linear--quadratic regulator
(LQR) problems with arbitrary stagewise linear equality constraints, using
linear work and logarithmic depth in the horizon for fixed local dimensions.
The solver uses a two-part reduction that separates feasibility from optimization.
First, an associative suffix scan computes the states from which each
remaining horizon is feasible.  Second, we use these feasible-continuation sets to
eliminate the constraints independently at each stage, producing an
unconstrained LQR problem with smaller or equal, possibly stage-varying dimensions.
A parallel unconstrained-LQR solver completes the primal solution.
Finally, we recover valid multipliers from the original KKT stationarity
equations using a logarithmic-depth contraction--expansion procedure, even
when the multipliers are not unique.
Sequential and parallel implementations
are validated against dense KKT solutions and original residuals.
We release a C++/CUDA implementation with JAX bindings under the MIT license.
\end{abstract}


\section{Introduction}
\label{sec:introduction}

\IEEEPARstart{E}{quality}-constrained LQR problems arise naturally in model predictive control
(MPC) and in the optimal control of differential-algebraic equations (DAEs).
In robotics, equality constraints enforce the kinematics of prescribed
contacts in legged locomotion, while inequality constraints enforce collision
avoidance and actuator limits.
In interior-point methods for MPC, eliminating slack and dual variables
associated with inequality constraints from the Newton--KKT system yields an LQR subproblem
subject to any remaining linearized stagewise
equalities~\cite{ref-InteriorPointMPC,ref-ParallelNewtonMPC}.
Active-set methods for MPC instead form equality-constrained LQR subproblems
by treating the inequalities in the current working set as equalities.
Equality-constrained LQR solves are therefore a core subroutine of active-set
methods even when the original stagewise constraints are all inequalities,
and of interior-point methods when stagewise equalities are present.

For nonlinear MPC, sequential quadratic programming (SQP) forms quadratic
subproblems that can be solved by these interior-point or active-set methods.
A nonlinear interior-point method can instead be applied directly to the
nonlinear problem, with the LQR subproblems arising from its primal--dual
Newton--KKT systems.

For semi-explicit DAEs, the algebraic variables can be included among the
controls.  Discretizing and linearizing the algebraic equations then produces
stagewise linear equalities of the form considered
here~\cite{ref-DAE-DynamicOptimization}.

For LQR problems with linear inequality constraints, projecting the feasible
set generally yields polyhedra rather than affine sets, so the bounded-size
representation used here does not carry over.
Such problems are typically solved by iterative methods, such as the
interior-point and active-set methods discussed above, with equality-constrained
LQR solves serving as subroutines.

Temporal sparsity makes these problems tractable.  Although the horizon-wide
optimality system is large, a Riccati recursion solves an unconstrained LQR
problem with work linear in the horizon.  Parallel formulations based on
associative scans reduce the horizon depth to logarithmic~\cite{ref-AS-LQR}.

Linear equality constraints complicate this picture.  A restriction at a
future state limits which present states admit a feasible continuation.
A sequential solver can propagate this information backward, eliminate the
controls that enforce it, and continue one stage at a time.  The resulting
work remains linear, but the propagation appears inherently sequential unless
the information carried by an interval is chosen carefully.

We separate feasibility from optimization without requiring linearly
independent equality constraints.  Specifically, our contributions are:
\begin{enumerate}
  \item A parallel scan that computes, for each stage, an affine
  parameterization of the states from which the remaining horizon is feasible.
  \item Independent stagewise constraint elimination producing an
  unconstrained LQR problem with smaller or equal state and control
  dimensions.  We prove a bijection between feasible trajectories of the
  original and reduced problems, objective equivalence, and preservation of
  standard convexity assumptions.
  \item A parallel method for recovering multipliers satisfying the original
  KKT equations, including when redundant constraints make the multipliers
  nonunique.
\end{enumerate}
For fixed local dimensions, the reduction and multiplier recovery require
linear work and logarithmic depth in the horizon.  Both also admit sequential
implementations with horizon-linear work.

\subsection{Related work}

Sideris and Rodriguez~\cite{ref-SR} treat both mixed and state-only equalities,
but their factorization is horizon-linear only when every control-side
constraint block has full row rank; state-only rows give cubic worst-case
horizon complexity.  Laine and Tomlin~\cite{ref-CLQR} give a horizon-linear
recursion for general equality constraints under standard LQR definiteness,
including constraint-aware feedback policies.

Vanroye, De Schutter, and Decré~\cite{ref-ConstrainedRiccati} recast the
recursion as a stagewise KKT factorization that avoids per-stage SVDs.  Its
positive-definite reduced-Hessian requirement admits indefinite exact
Lagrangian Hessians, and it underlies FATROP~\cite{ref-FATROP}.  It assumes a
full-row-rank global equality Jacobian; dynamically mediated redundancies
require global rather than stagewise rank-revealing preprocessing.

Yang et al.~\cite{ref-FactorGraphLQR} propagate cost-to-go and constraint-to-go
factors by backward factor-graph elimination, equivalently to constrained
dynamic programming for stage-local constraints.

Feasible-continuation sets also arise in constrained dynamic
programming~\cite[Sec.~2.3]{ref-MPC-Book}.  Here, associative composition of
bounded-size affine relations enables their work-efficient parallel
construction.

Lefebvre and Vantilborgh~\cite{ref-ParallelEC-OCP} use parallel-prefix
conditional-value tuples propagating feasibility and optimality together.
Our cost-independent feasibility scan and stagewise elimination instead
produce an unconstrained LQR problem whose parallel Riccati passes operate on
reduced state and control dimensions.

\section{Problem formulation}
\label{sec:problem}

For stages $i=0,\ldots,N-1$, let $x_i\in\R^{n_i}$ denote the state and
$u_i\in\R^{m_i}$ the control, and let $x_N\in\R^{n_N}$ denote the terminal
state.  Define the stage cost
\begin{equation}
\label{eq:stage-cost}
\begin{aligned}
\ell_i(x_i,u_i)={}&\tfrac12x_i^\trans Q_ix_i+x_i^\trans M_iu_i
+\tfrac12u_i^\trans R_iu_i\\
&+q_i^\trans x_i+r_i^\trans u_i.
\end{aligned}
\end{equation}
We consider
\begin{equation}
\label{eq:problem}
\begin{aligned}
\min_{x,u}\quad
&\sum_{i=0}^{N-1}\ell_i(x_i,u_i)+\tfrac12x_N^\trans Q_Nx_N+q_N^\trans x_N\\
\text{s.t.}\quad
&\begin{aligned}[t]
&x_0=s_0,\\
&x_{i+1}=A_ix_i+B_iu_i+c_i,
&&i=0,\ldots,N-1,\\
&C_ix_i+D_iu_i+d_i=0,
&&i=0,\ldots,N-1,\\
&E_ix_i+e_i=0,
&&i=0,\ldots,N.
\end{aligned}
\end{aligned}
\end{equation}
The minimization variables $x=(x_0,\ldots,x_N)$ and
$u=(u_0,\ldots,u_{N-1})$ collect the trajectory.
Empty matrices represent absent constraints.  Dimensions and constraint ranks
may vary by stage.  We do not assume that the rows of $[C_i\ D_i]$ or $E_i$
are independent.

The standard convex LQR assumptions are
\begin{equation}
\label{eq:convexity}
\begin{bmatrix}Q_i&M_i\\M_i^\trans&R_i\end{bmatrix}\succeq0,
\qquad R_i\succ0,
\qquad Q_N\succeq0.
\end{equation}
Here $\succ0$ and $\succeq0$ denote positive definiteness and positive
semidefiniteness, respectively.  The matrices $Q_i$ and $R_i$ are symmetric.
The reduction does not require \eqref{eq:convexity} and remains
valid for arbitrary symmetric quadratic costs.  When \eqref{eq:problem} is
feasible, \eqref{eq:convexity} is a sufficient, but not necessary, condition
for the existence of a unique minimizer.

\section{Affine relations}
\label{sec:relations}

\begin{definition}[Relation]
\label{def:relation}
A relation between two sets is a subset of their Cartesian product.
\end{definition}

\begin{definition}[Relational composition]
\label{def:relational-composition}
Let $\mathcal R$ and $\mathcal S$ be relations on
$\R^p\times\R^q$ and $\R^q\times\R^r$, respectively.  Their relational
composition is defined as
\begin{equation}
\label{eq:composition-set}
\mathcal R\circ\mathcal S
:=\{(\xi,\zeta):\exists\eta,\ (\xi,\eta)\in\mathcal R\ \wedge\ (\eta,\zeta)\in\mathcal S\}.
\end{equation}
\end{definition}

\begin{proposition}[Associativity]
\label{prop:associativity}
For relations with compatible endpoint spaces, composition is associative:
\begin{equation}
(\mathcal R\circ\mathcal S)\circ\mathcal T
=\mathcal R\circ(\mathcal S\circ\mathcal T).
\end{equation}
\end{proposition}

\begin{proof}
For any endpoint pair $(\xi,\omega)$, membership in either side is equivalent
to the existence of $\eta,\zeta$ such that $(\xi,\eta)\in\mathcal R$,
$(\eta,\zeta)\in\mathcal S$, and $(\zeta,\omega)\in\mathcal T$.
\end{proof}

\begin{definition}[Affine relation]
An affine relation is a set
\begin{equation}
\label{eq:relation}
\mathcal R
:=\{(\xi,\eta)\in\R^p\times\R^q:U\xi+V\eta+w=0\},
\end{equation}
where $U$, $V$, and $w$ are fixed and have compatible dimensions.
Different triples $(U,V,w)$ that describe the same set are regarded as
equivalent.
\end{definition}

\begin{lemma}[Existential projection]
\label{lem:projection}
Let $F y+G z+g=0$, and let the columns of $W$ form a basis of
$\nullsp(F^\trans)$.  Then
\begin{equation}
\label{eq:projection}
\exists y:\;Fy+Gz+g=0
\quad\Longleftrightarrow\quad
W^\trans(Gz+g)=0.
\end{equation}
If a row reduction produces $0=\gamma$ with $\gamma\ne0$, the projected set is
empty.
\end{lemma}

\begin{proof}
The equation in $y$ is solvable exactly when $-(Gz+g)\in\range(F)$.  By the
fundamental theorem of linear algebra,
$\range(F)=\nullsp(F^\trans)^\perp$, which is equivalent to
\eqref{eq:projection}.
\end{proof}

\begin{proposition}[Closure]
\label{prop:closure}
The composition of two affine relations with compatible endpoint spaces
is an affine relation.
\end{proposition}

\begin{proof}
Writing $\mathcal R$ and $\mathcal S$ with coefficient triples
$(U_1,V_1,w_1)$ and $(U_2,V_2,w_2)$, respectively, apply
\cref{lem:projection} to eliminate $\eta$ from
\begin{equation}
\label{eq:composition-matrix}
\underbrace{\begin{bmatrix}V_1\\U_2\end{bmatrix}}_{F}\eta
+\underbrace{\begin{bmatrix}U_1&0\\0&V_2\end{bmatrix}}_{G}
\begin{bmatrix}\xi\\\zeta\end{bmatrix}
+\underbrace{\begin{bmatrix}w_1\\w_2\end{bmatrix}}_{g}=0.\qedhere
\end{equation}
\end{proof}

\section{Computing feasible-continuation sets}
\label{sec:domains}

Below, we use a dummy terminal endpoint $x_{N+1}\in\R^{n_{N+1}}$,
where $n_{N+1}=0$.

\subsection{Interval feasibility relations}

\begin{definition}[Interval feasibility relation]
\label{def:interval-relation}
For $0\le i<j\le N+1$,
\[
\begin{aligned}
&\mathcal R_{i\to j}:=\bigl\{(x_i,x_j)\in\R^{n_i}\times\R^{n_j}\bigm|\\
&\quad\exists\,(x_k)_{i<k<j},\ (u_k)_{i\le k<\min(j,N)}\ \text{such that}\\
&\bigl(\forall k\in\{i,\ldots,\min(j,N)-1\},\;
x_{k+1}=A_kx_k+B_ku_k+c_k\bigr)\\
&\wedge\bigl(\forall k\in\{i,\ldots,\min(j,N)-1\},\;
0=C_kx_k+D_ku_k+d_k\bigr)\\
&\wedge\bigl(\forall k\in\{i,\ldots,j-1\},\;
0=E_kx_k+e_k\bigr)\bigr\}.
\end{aligned}
\]
\end{definition}

\begin{theorem}[Interval composition]
\label{thm:interval-composition}
For $0\le i<j\le N+1$,
\begin{equation}
\label{eq:interval-composition}
\mathcal R_{i\to j}
=\mathcal R_{i\to i+1}\circ\cdots\circ\mathcal R_{j-1\to j}.
\end{equation}
\end{theorem}

\begin{proof}
The case $j=i+1$ is trivial. For $j>i+1$, the constraints in
\cref{def:interval-relation} at stage $i$ share only $x_{i+1}$ with those
at stages $i+1,\ldots,j-1$. By \eqref{eq:composition-set},
$\mathcal R_{i\to j}=\mathcal R_{i\to i+1}\circ\mathcal R_{i+1\to j}$.
Induction on $j-i$ and \cref{prop:associativity} prove
\eqref{eq:interval-composition}.
\end{proof}

\begin{corollary}[Affine interval feasibility relations]
\label{cor:affine-interval-relations}
For every $0\le i<j\le N+1$, $\mathcal R_{i\to j}$ is an affine relation.
\end{corollary}

\begin{proof}
For $i<N$, stack the stage dynamics and constraints as
\begin{equation}
\label{eq:base-system}
\underbrace{\begin{bmatrix}B_i\\D_i\\0\end{bmatrix}}_{F_i}u_i
+\underbrace{\begin{bmatrix}
A_i&-I\\C_i&0\\E_i&0
\end{bmatrix}}_{G_i}
\begin{bmatrix}x_i\\x_{i+1}\end{bmatrix}
+\underbrace{\begin{bmatrix}c_i\\d_i\\e_i\end{bmatrix}}_{g_i}=0.
\end{equation}
The single-stage feasibility relation is therefore
\begin{equation}
\label{eq:base-relation}
\mathcal R_{i\to i+1}
=\biggl\{(x_i,x_{i+1}):\exists u_i,\ F_i u_i+G_i
\begin{bmatrix}x_i\\x_{i+1}\end{bmatrix}+g_i=0\biggr\}.
\end{equation}
By \cref{lem:projection}, eliminating $u_i$ gives an affine representation.
The terminal relation
$\mathcal R_{N\to N+1}=\{(x_N,x_{N+1}):E_Nx_N+e_N=0\}$ is also affine.
The result follows from \cref{prop:closure,thm:interval-composition}.
\end{proof}

Each base relation depends only on its stage's data, so all can be
constructed in parallel.

\subsection{Feasible-continuation sets}

\begin{definition}[Feasible-continuation set]
\label{def:feasible-continuation-set}
The feasible-continuation set $\mathcal X_i$ consists of the states $x_i$
with $(x_i,x_{N+1})\in\mathcal R_{i\to N+1}$.
\end{definition}

By \cref{cor:affine-interval-relations}, $\mathcal X_i$ has an affine representation
\begin{equation}
\label{eq:domain-equations}
\mathcal X_i=\{x_i:H_ix_i+h_i=0\}.
\end{equation}

By \cref{prop:associativity}, all $N+1$ suffix compositions in
\eqref{eq:interval-composition} with $j=N+1$ can be computed by a work-efficient associative
suffix scan~\cite{ref-AS}.  A balanced up-sweep composes interval relations;
the down-sweep distributes the appropriate right aggregate to each interval.
There are $O(N)$ relation compositions and $O(\log N)$ composition levels.

Each nonempty $\mathcal X_i$ admits a compact representation with at most
$n_i$ rows in $H_i$, obtained by retaining only independent equations.

\subsection{Feasible-continuation set parameterization}

If any $\mathcal X_i$ is represented by an inconsistent affine equation, no
state at that stage admits a feasible continuation.

Otherwise, choose $T_i$ and $t_i$ such that
\begin{equation}
\label{eq:state-param}
H_iT_i=0,
\qquad H_it_i+h_i=0,
\qquad \rank(T_i)=\bar n_i,
\end{equation}
where the columns of $T_i\in\R^{n_i\times\bar n_i}$ form a basis of
$\nullsp(H_i)$.  Thus $\mathcal X_i=\{T_iz_i+t_i:z_i\in\R^{\bar n_i}\}$.
Choose a left inverse $S_iT_i=I_{\bar n_i}$, giving
$z_i=S_i(x_i-t_i)$ on $\mathcal X_i$.  These factorizations are independent
across stages after the scan.

The original problem is feasible if and only if $H_0s_0+h_0=0$.
Its reduced initial state is the unique coordinate $z_0=S_0(s_0-t_0)$.

\section{Stage-parallel constraint elimination}
\label{sec:local-elimination}

All parameterizations $(T_i,t_i,S_i)$ are now known.  The construction below
uses only original stage $i$ and the feasible-continuation sets $\mathcal X_i$ and
$\mathcal X_{i+1}$; hence all stages can be processed simultaneously.

\subsection{Feasible controls}

Substitute $x_i=T_iz_i+t_i$.  A control must satisfy the original mixed
constraint and must make the dynamics successor belong to
$\mathcal X_{i+1}$.  Using \eqref{eq:domain-equations}, these requirements are
\begin{equation}
\label{eq:local-control-system}
J_i u_i+K_i z_i+k_i=0,
\end{equation}
where
\begin{equation}
\label{eq:local-control-matrices}
\begin{aligned}
J_i&=
\begin{bmatrix}D_i\\H_{i+1}B_i\end{bmatrix},&
K_i&=
\begin{bmatrix}C_iT_i\\H_{i+1}A_iT_i\end{bmatrix},\\
k_i&=
\begin{bmatrix}
C_it_i+d_i\\H_{i+1}(A_it_i+c_i)+h_{i+1}
\end{bmatrix}.
\end{aligned}
\end{equation}
The state-only equality $E_ix_i+e_i=0$ need not be repeated: it is already
encoded in $\mathcal X_i$.

\begin{lemma}[Stagewise control parameterization]
\label{lem:control-param}
For every $z_i$, there exists $u_i$ satisfying
\eqref{eq:local-control-system}.  Moreover,
there exist $Y_i$, $Z_i$, and $y_i$ such that all its solutions are
\begin{equation}
\label{eq:control-param}
u_i=Y_iz_i+Z_iv_i+y_i,
\qquad v_i\in\R^{\bar m_i},
\end{equation}
where the columns of $Z_i$ form a basis of $\nullsp(J_i)$.
\end{lemma}

\begin{proof}
For any $z_i$, $T_iz_i+t_i\in\mathcal X_i$.  By
\cref{def:interval-relation}, some control satisfies the stage mixed constraint
and sends the state into $\mathcal X_{i+1}$; this is exactly
\eqref{eq:local-control-system}.  Consequently,
$\range(K_i)\subseteq\range(J_i)$ and
$k_i\in\range(J_i)$.  Hence there exist $Y_i,y_i$ satisfying
$J_iY_i=-K_i$ and $J_i y_i=-k_i$.
Thus $u_i=Y_iz_i+y_i$ satisfies
\eqref{eq:local-control-system}.  Every other solution differs from it by an
element of $\nullsp(J_i)$, giving \eqref{eq:control-param}.
\end{proof}

\subsection{Reduced dynamics}

The dynamics successor produced by \eqref{eq:control-param} lies in
$\mathcal X_{i+1}$.  Its unique reduced coordinate is obtained with
$S_{i+1}$:
\begin{equation}
\label{eq:reduced-dynamics}
z_{i+1}=\bar A_i z_i+\bar B_i v_i+\bar c_i,
\end{equation}
where
\begin{equation}
\label{eq:reduced-dynamics-data}
\begin{aligned}
\bar A_i&=S_{i+1}(A_iT_i+B_iY_i),\\
\bar B_i&=S_{i+1}B_iZ_i,\\
\bar c_i&=S_{i+1}(A_it_i+B_iy_i+c_i-t_{i+1}).
\end{aligned}
\end{equation}

\subsection{Reduced objective}

First substitute the state parameterization and define
\begin{equation}
\label{eq:state-cost-data}
\begin{aligned}
Q_i^z&=T_i^\trans Q_iT_i,&
M_i^z&=T_i^\trans M_i,\\
q_i^z&=T_i^\trans(Q_it_i+q_i),&
r_i^z&=M_i^\trans t_i+r_i.
\end{aligned}
\end{equation}
\begingroup
% Keep the display with its lead-in, with slightly tighter inter-row spacing.
\setlength{\jot}{2pt}
Substitution of \eqref{eq:control-param} then gives a reduced stage cost of the
standard LQR form, with
\begin{equation}
\label{eq:reduced-cost-data}
\begin{aligned}
\bar Q_i={}&Q_i^z+M_i^zY_i+Y_i^\trans(M_i^z)^\trans
             +Y_i^\trans R_iY_i,\\
\bar M_i={}&M_i^zZ_i+Y_i^\trans R_iZ_i,\\
\bar R_i={}&Z_i^\trans R_iZ_i,\\
\bar q_i={}&q_i^z+M_i^zy_i+Y_i^\trans(R_iy_i+r_i^z),\\
\bar r_i={}&Z_i^\trans(R_iy_i+r_i^z).
\end{aligned}
\end{equation}
\endgroup
The discarded constant is
\begin{equation}
\label{eq:stage-constant}
\kappa_i=\tfrac12t_i^\trans Q_it_i+q_i^\trans t_i
+\tfrac12y_i^\trans R_iy_i+(r_i^z)^\trans y_i.
\end{equation}
At the terminal stage,
\begin{equation}
\label{eq:terminal-cost-data}
\bar Q_N=T_N^\trans Q_NT_N,
\qquad
\bar q_N=T_N^\trans(Q_Nt_N+q_N),
\end{equation}
with constant
$\kappa_N=\tfrac12t_N^\trans Q_Nt_N+q_N^\trans t_N$.

Together with $z_0=S_0(s_0-t_0)$,
\labelcref{eq:reduced-dynamics,eq:reduced-cost-data,eq:terminal-cost-data}
define an unconstrained LQR problem in $(z_i,v_i)$.

\begin{theorem}[Equivalence]
\label{thm:equivalence}
Assume $s_0\in\mathcal X_0$.  The maps
\begin{equation}
\label{eq:recovery-maps}
x_i=T_iz_i+t_i,
\qquad
u_i=Y_iz_i+Z_iv_i+y_i
\end{equation}
give a bijection between trajectories of the reduced LQR problem and feasible
trajectories of \eqref{eq:problem}.  Corresponding objective values differ by
the constant $\sum_{i=0}^N\kappa_i$.  Hence corresponding minimizers coincide.
\end{theorem}

\begin{proof}
For a reduced trajectory, the state parameterization enforces every state-only
constraint.  Equation \eqref{eq:local-control-system} enforces the original
mixed constraint and places the dynamics successor in
$\mathcal X_{i+1}$.  Since $S_{i+1}$ is a left inverse of $T_{i+1}$,
\eqref{eq:reduced-dynamics} makes the coordinate of this successor equal to
the declared $z_{i+1}$; therefore the original dynamics hold.

Conversely, every feasible state has a unique coordinate $z_i$ because $T_i$
has full column rank.  For this coordinate, every feasible control belongs to
the complete solution set \eqref{eq:control-param}, and $Z_i$ has full column
rank, so it has a unique coordinate $v_i$.  The original dynamics then imply
\eqref{eq:reduced-dynamics}.  Direct substitution gives
\eqref{eq:reduced-cost-data} and the stated constant.
\end{proof}

\begin{proposition}[Preservation of convexity]
\label{prop:convexity}
Under \eqref{eq:convexity},
the reduced stage Hessians and $\bar Q_N$ are positive semidefinite,
and $\bar R_i\succ0$ whenever $\bar m_i>0$.
\end{proposition}

\begin{proof}
The reduced stage Hessian is
\[
\begin{bmatrix}T_i&0\\Y_i&Z_i\end{bmatrix}^\trans
\begin{bmatrix}Q_i&M_i\\M_i^\trans&R_i\end{bmatrix}
\begin{bmatrix}T_i&0\\Y_i&Z_i\end{bmatrix}\succeq0.
\]
Furthermore,
$\bar R_i=Z_i^\trans R_iZ_i\succ0$ because $R_i\succ0$ and $Z_i$ has full
column rank.  The terminal statement follows by congruence with $T_N$.
\end{proof}

\section{Multiplier recovery}
\label{sec:multipliers}

The primal reduction does not need to retain the row operations associated
with the original constraints.  Instead, multipliers can be recovered directly
from the original KKT stationarity equations after $(x,u)$ is known.

Use multipliers $\lambda_{i+1}$ for the dynamics residual
$A_ix_i+B_iu_i+c_i-x_{i+1}=0$, $\mu_i$ for the mixed constraint,
$\eta_i$ for the state-only constraint, and $\lambda_0$ for the initial
constraint $s_0-x_0=0$.  Define
\begin{equation}
\label{eq:objective-gradients}
\begin{aligned}
g_i^x&=Q_ix_i+M_iu_i+q_i,&
g_i^u&=M_i^\trans x_i+R_iu_i+r_i,\\
g_N^x&=Q_Nx_N+q_N.
\end{aligned}
\end{equation}
State and control stationarity at stage $i$ become
\begin{equation}
\label{eq:stage-stationarity}
\begin{aligned}
g_i^x-\lambda_i+A_i^\trans\lambda_{i+1}
+C_i^\trans\mu_i+E_i^\trans\eta_i&=0,\\
g_i^u+B_i^\trans\lambda_{i+1}+D_i^\trans\mu_i&=0,
\end{aligned}
\end{equation}
and terminal stationarity is
\begin{equation}
\label{eq:terminal-stationarity}
g_N^x-\lambda_N+E_N^\trans\eta_N=0.
\end{equation}

Equations \eqref{eq:stage-stationarity} define an affine relation between
$\lambda_i$ and $\lambda_{i+1}$.  Indeed, stack them as
\begin{equation}
\label{eq:dual-base-system}
\underbrace{\begin{bmatrix}
C_i^\trans&E_i^\trans\\D_i^\trans&0
\end{bmatrix}}_{F_i^\lambda}
\begin{bmatrix}\mu_i\\\eta_i\end{bmatrix}
+\underbrace{\begin{bmatrix}
-I&A_i^\trans\\0&B_i^\trans
\end{bmatrix}}_{G_i^\lambda}
\begin{bmatrix}\lambda_i\\\lambda_{i+1}\end{bmatrix}
+\underbrace{\begin{bmatrix}g_i^x\\g_i^u\end{bmatrix}}_{g_i^\lambda}=0.
\end{equation}
Eliminating $(\mu_i,\eta_i)$ by \cref{lem:projection} gives a dual endpoint
relation $\mathcal D_{i\to i+1}$ on $(\lambda_i,\lambda_{i+1})$.
Likewise, eliminating $\eta_N$ from
\eqref{eq:terminal-stationarity} gives a terminal relation
$\mathcal D_{N\to N+1}$ on $(\lambda_N,\lambda_{N+1})$, with a
zero-dimensional dummy endpoint $\lambda_{N+1}\in\R^0$.

The dual relations can be contracted by the same balanced composition used in
\cref{sec:domains}.  To compose $\mathcal D_{i\to k}$ and
$\mathcal D_{k\to j}$, with $i<k<j$, stack their equations as in
\eqref{eq:composition-matrix} and eliminate $\lambda_k$.  The residual equations
define $\mathcal D_{i\to j}$.  Solving the stacked equations for the eliminated
costate gives an affine function of the endpoint costates and the remaining
free variables, which we collect in $\rho_{i\to j}$.  Denoting its coefficients
by $\Psi_{i\to j}$, $\psi_{i\to j}$, and $\Omega_{i\to j}$, we obtain and
retain the recovery map
\begin{equation}
\label{eq:middle-dual-map}
\lambda_k=\Psi_{i\to j}
\begin{bmatrix}\lambda_i\\\lambda_j\end{bmatrix}
+\psi_{i\to j}+\Omega_{i\to j}\rho_{i\to j}.
\end{equation}
For every $(\lambda_i,\lambda_j)\in\mathcal D_{i\to j}$, this parameterizes
all compatible shared costates; the columns of $\Omega_{i\to j}$ form a
basis of the nullspace of the coefficient matrix of $\lambda_k$ in the
stacked system.

After contraction, choose any $\lambda_0$ in the root feasible set; a
minimum-norm particular solution is a natural choice.  Traverse the stored
tree from root to leaves, setting the free coordinates $\rho_{i\to j}$ to zero.  Each
tree level is independent, so all costates are recovered in $O(\log N)$
horizon depth.  Given adjacent costates, solve
\eqref{eq:dual-base-system} independently at every stage for
$(\mu_i,\eta_i)$, and solve \eqref{eq:terminal-stationarity} for $\eta_N$.

\begin{samepage}
\begin{proposition}[Correctness of multiplier recovery]
\label{prop:multipliers}
Let $(x,u)$ be the trajectory recovered from an optimum of the reduced LQR
problem.
The contraction--expansion procedure above returns multipliers satisfying all
KKT equations of \eqref{eq:problem}.
\end{proposition}
\end{samepage}

\begin{proof}
By \cref{thm:equivalence}, $(x,u)$ is feasible and optimal.  Its objective
gradient is orthogonal to every feasible direction, hence belongs to the
range of the stacked constraint transpose; KKT multipliers therefore exist
without convexity.  Each dual base relation results from exact elimination
of the local multipliers from the stationarity equations.
Contraction preserves its feasible endpoints, expansion
selects compatible internal costates, and the leaf solves satisfy the
original equations.  All original stationarity equations consequently hold.
\end{proof}

If the original equality rows are redundant, multipliers need not be unique.
The free variables in \eqref{eq:middle-dual-map} and the leaf systems expose
this nonuniqueness explicitly; setting them to zero selects one representative.
For a sequential implementation, a representative can be obtained more
simply by propagating affine sets of admissible costates right-to-left and then
substituting forward.  That variant has $O(N)$ horizon depth but uses the same
local equations.

\section{Complete algorithm and complexity}
\label{sec:algorithm}

\Cref{alg:complete} summarizes the complete solver, from the feasibility scan
through primal reconstruction and multiplier recovery.

\begin{algorithm}[!b]
\caption{Parallel equality-constrained LQR solver}
\label{alg:complete}
\begin{algorithmic}
\REQUIRE Data of \eqref{eq:problem}
\STATE Build base relations \eqref{eq:base-relation} in parallel.
\STATE Compute all suffix relations by an associative scan.
\STATE Return infeasible if a suffix is inconsistent or
$H_0s_0+h_0\ne0$.
\STATE In parallel, construct $(T_i,t_i,S_i)$.
\STATE In parallel, form \eqref{eq:control-param} and the reduced
dynamics and costs (\cref{sec:local-elimination}).
\STATE Solve the reduced LQR problem by a parallel scan.
\STATE Recover $(x_i,u_i)$ in parallel by \eqref{eq:recovery-maps}.
\STATE Recover multipliers by contraction--expansion
(\cref{sec:multipliers}).
\end{algorithmic}
\end{algorithm}

Let $n=\max_i n_i$, $m=\max_i m_i$, and let $p$ bound the
total number of mixed and state-only equality rows at a stage.
A nonempty affine relation between interval endpoints has at most $2n$
independent rows.
Factoring a base system costs $O((n+p)(n+m)^2)$; an interval
composition eliminates $n$ shared variables from at most $4n$ canonical rows,
costing $O(n^3)$.  Local elimination and cost transformation cost
$O((n+p)(n+m)^2+(n+m)^3)$ per stage with dense kernels, giving the bound
$O(N[(n+p)(n+m)^2+(n+m)^3])$ for the primal solution.  When $p=O(n+m)$,
local work is cubic in $n+m$.
For fixed local dimensions, storage is $O(N)$ and the scan has $O(\log N)$
horizon depth; state parameterization and stage elimination have constant
horizon depth.  Local factorization depth multiplies the scan depth when
counting scalar operations.

The reduced Riccati recursion exploits stage-dependent dimensions
$(\bar n_i,\bar m_i)$.  It has linear horizon work and depth sequentially;
a work-efficient parallel unconstrained-LQR scan~\cite{ref-AS-LQR} reduces
the depth to $O(\log N)$.

For multiplier recovery, the local systems \eqref{eq:dual-base-system}
have at most $n+m$ equations and $p$ local multiplier variables.
Building the dual base relations and recovering these multipliers costs
$O((n+p)(n+m)^2)$ per stage.  Forming each dual composition and its recovery
map costs $O(n^3)$; evaluating the map during expansion costs $O(n^2)$.
The balanced contraction--expansion tree has $O(N)$ nodes and $O(\log N)$
levels, and the local multiplier systems are solved independently across
stages.  Thus the complete solver, including multiplier recovery, retains
the stated work bound and, for fixed local dimensions, $O(N)$ storage and
$O(\log N)$ horizon depth.

There is a useful shortcut.  If all state-only equalities are vacuous, the
terminal equality is vacuous, and eliminating $u_i$ from each mixed equality
leaves no restriction on $x_i$, then every
$\mathcal X_i=\R^{n_i}$.  The scan in \cref{sec:domains} can then be skipped,
and the mixed constraints can be eliminated locally.

\section{Implementation}
\label{sec:implementation}

Our MIT-licensed C++/CUDA package~\cite{ref-CLQR-Elimination}, with a native
JAX interface, implements sequential and parallel variants.
The native CPU path fuses suffix propagation and elimination;
CUDA uses compact, active-dimension kernels.  Intermediate dimensions depend
on data-dependent constraint ranks.  Under JAX's static-shape compilation,
a pure-JAX implementation must pad these intermediates to fixed upper bounds,
performing dense algebra in larger dimensions than the reduction requires.
We therefore implement the GPU kernels in CUDA and expose them through JAX,
retaining active-dimension computation within compiled JAX workflows.

Tests cover mixed, state-only, rank-deficient, and redundant constraints.
They check agreement with dense KKT solutions, original KKT residuals, and
objective equivalence (including constant terms).

\section{Benchmarks}
\label{sec:benchmarks}

Let $p_s$ and $p_m$ denote the numbers of state-only and mixed
state/control equality rows per stage, respectively.

% Use the plot markers at a uniform size, not font-dependent math symbols.
\newcommand{\benchmarkMarker}[1]{%
\raisebox{-1bp}{\includegraphics[page=#1,width=6bp]{figures/benchmark_markers.pdf}}}
% Keep citations in LaTeX so legend links and numbering follow the bibliography.
\newcommand{\benchmarkLegend}[1]{%
{\fontsize{7}{8}\selectfont\makebox[\columnwidth][c]{%
\benchmarkMarker{1}\,Ours (CPU)\quad
\benchmarkMarker{2}\,\cite{ref-ConstrainedRiccati}\quad
\benchmarkMarker{3}\,\cite{ref-FactorGraphLQR}\quad
\benchmarkMarker{4}\,\cite{ref-CLQR}\quad
\benchmarkMarker{5}\,#1}}\par\vspace{3pt}}

\begin{figure}[!hb]
\centering
\benchmarkLegend{Ours (CUDA/JAX)}
\includegraphics[width=0.95\columnwidth]{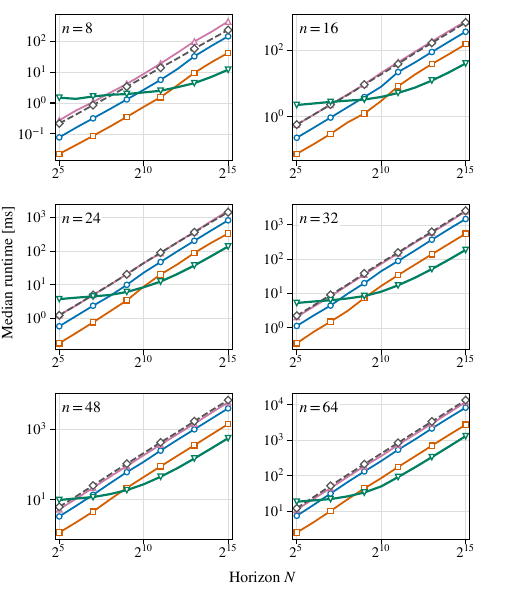}
\caption{Median solve time versus horizon at fixed state dimension.}
\label{fig:benchmark-runtime}
\medskip
\benchmarkLegend{Ours (CUDA/JAX)}
\includegraphics[width=0.95\columnwidth]{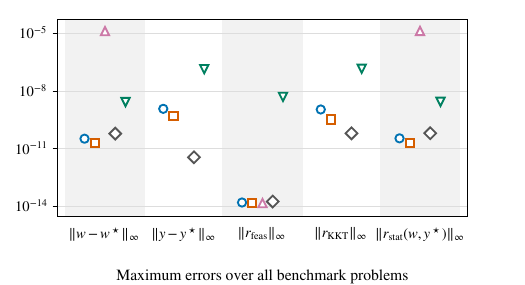}
\caption{Maximum errors over all benchmark problems.
$w=(x,u)$ and $y$ collect the primal variables and equality multipliers;
$(w^\star,y^\star)$ is the known optimum.
The factor-graph method~\cite{ref-FactorGraphLQR} returns no multipliers,
so its dual error and KKT residual are unavailable.
The last group independently audits stationarity using $y^\star$ for every method.}
\label{fig:benchmark-accuracy}
\end{figure}

The FP64 measurements in \cref{fig:benchmark-runtime,fig:benchmark-accuracy}
used an Ubuntu 24.04.4 workstation with a 24-core Intel Core Ultra 9 285K CPU
and an NVIDIA RTX PRO 6000 Blackwell Workstation Edition GPU
(95.6 GiB reported device memory), using GCC 13.3, CUDA 12.9, and NVIDIA
driver 610.43.02.  The CPU solvers use native-architecture Release
optimization and one thread.
All solvers receive identical dense instances constructed with known
optimal states, controls, and multipliers, with
$m=n/2$, $p_s=n/4$ state-only rows, and $p_m=n/8$ mixed rows.
The comparison focuses on scaling with the horizon at fixed state dimension.
We use every combination of $N=2^5,2^6,\ldots,2^{15}$ and
$n\in\{8,16,24,32,48,64\}$.
State equalities apply at nodes $1,\ldots,N$, avoiding redundancy at the
fixed initial state.

We use the authors' implementations of generalized Riccati
recursion~\cite{ref-ConstrainedRiccati}, factor-graph
elimination~\cite{ref-FactorGraphLQR}, and constrained dynamic
programming~\cite{ref-CLQR}.
Our CPU solver uses native dense kernels, without an external BLAS dependency.
Our solvers and the implementations of~\cite{ref-ConstrainedRiccati,ref-CLQR}
recover multipliers inside the timed solve; the implementation
of~\cite{ref-FactorGraphLQR} does not return multipliers.

Each timed measurement includes one factorization and one solve.
Runtimes are medians of
timed solves accumulated for at least one second, after one untimed warm-up solve.
Preparation of input representations and one-time workspace allocation are
excluded; the script also reports setup-plus-solve timings.
JAX uses device-resident numerical inputs and outputs, excludes compilation,
and synchronizes each repetition.  Generation and validation are untimed.

The unscaled original KKT residual covers the initial condition, dynamics,
equality feasibility, and stationarity with returned multipliers.
Primal and dual errors use the known optimum.  For primal-only solvers,
stationarity is separately audited with the known optimal multipliers;
this is not a solver-produced KKT residual.

The method of~\cite{ref-ConstrainedRiccati} is fastest on CPU within its
full-row-rank regime; our method also supports redundant global equalities
and nonunique multipliers.
Our sequential solver is faster than the measured implementations
of~\cite{ref-CLQR,ref-FactorGraphLQR} on all benchmarked cases.
At $N=32768$, $n=8$, the device-resident GPU solve takes $11.965$ ms, versus
$41.932$ ms for the method
of~\cite{ref-ConstrainedRiccati}, a $3.50\times$ speedup over the fastest
measured CPU method.  At this horizon, the speedup ranges from $2.09\times$
to $3.82\times$ across all tested state dimensions.
For longer horizons, the fixed GPU's finite concurrency and linear data
volume eventually dominate the logarithmic scan depth; neither constant
kernel time nor unbounded speedup is implied.

\section{Conclusion}

We have presented a parallel solver for linear-equality-constrained LQR
problems that separates feasibility from optimization.  An associative scan
computes the states from which each remaining horizon is feasible.  These
sets allow the stagewise constraints to be eliminated independently,
yielding an unconstrained LQR problem with smaller or equal state and control
dimensions.  The construction accommodates mixed and state-only constraints,
including redundant equalities.

The same affine-relation framework recovers multipliers from the original
KKT stationarity equations, including when they are nonunique.  The complete
solver therefore obtains both primal and dual solutions with linear work
and logarithmic depth in the horizon for fixed local dimensions.

The benchmarks demonstrate benefits for both sequential and parallel
implementations.  Our sequential solver is faster than the implementations
of~\cite{ref-CLQR,ref-FactorGraphLQR} on all benchmarked cases.  At the
longest tested horizon, the device-resident GPU solver outperforms the
fastest measured CPU method across all tested state dimensions.
The C++/CUDA implementation and its JAX bindings are available under the
MIT License~\cite{ref-CLQR-Elimination}.

\bibliographystyle{ieeetr}
\bibliography{references}

\end{document}